\documentclass[12pt]{article}
\usepackage{amsfonts}
\usepackage{amssymb}
\usepackage{amsmath}
\usepackage{epsfig}
\usepackage{amsthm}

\newtheorem{defi}{Definition}[section]
\newtheorem{proposition}[defi]{Proposition}
\newtheorem{theorem}[defi]{Theorem}
\newtheorem{lemma}[defi]{Lemma}
\newtheorem{exa}[defi]{Example}
\newtheorem{corollary}[defi]{Corollary}
\newtheorem{remark}[defi]{Remark}
\newcommand{\Gal}{\mbox{\rm Gal}}

\newcommand{\Aut}{\mbox{\rm Aut}}

\newcommand{\llll}{\lambda}
\newcommand{\circo}{{\#}}
\newcommand{\Z}{{\mathbb Z}}
\newcommand{\Q}{{\mathbb Q}}
\newcommand{\Qpt}{{\mathbb Q}_p^{\rm \scriptsize tame}}
\newcommand{\Qpur}{{\mathbb Q}_p^{\rm \scriptsize ur}}
\newcommand{\ex}{{\rm \scriptsize ur}}
\newcommand{\sdp}{\rtimes}

\newcommand{\Pk}{{\mathbb P}^{1}_{k}}
\newcommand{\la}{\longrightarrow}

\newcommand{\PK}{{\mathbb P}^{1}_{K}}
\newcommand{\PKK}{{\mathbb P}^{1}_{K'}}

\title{{Wild monodromy and automorphisms of curves}
\thanks{Work supported in part by the European Community's Sixth Framework
Programme under Contract MRTN-CT2003-504917}
}
\author{Claus Lehr and Michel Matignon}
\begin{document} 
\maketitle
\begin{abstract} 

Let $R$ be a complete discrete valuation ring of mixed characteristic $(0,p)$
with field of fractions $K$ containing the $p$-th roots of unity.
This paper is concerned with semi-stable models of $p$-cyclic covers
of the projective line $C \la \PK$.
We start by providing a new construction of a semi-stable model 
of $C$ in the case of an equidistant branch locus.
If the cover is given by the Kummer equation $Z^p=f(X_0)$ we 
define what we called the monodromy polynomial
${\mathcal L}(Y)$ of $f(X_0)$;
a polynomial with coefficients in $K$.
Its zeros are key to obtaining a semi-stable model of $C$.
As a corollary we obtain an upper bound for the minimal 
extension $K'/K$ over which a stable model of the curve $C$ exists.
Consider the polynomial ${\cal L}(Y)\prod(Y^p-f(y_i))$ where the $y_i$ 
range over the zeros of ${\cal L}(Y)$.
We show that the splitting field of this polynomial always contains
$K'$, and that in some instances the two fields are equal.
\end{abstract}

\section{Introduction}

Let $R$ be a complete discrete valuation ring of mixed characteristic $(0,p)$
with field of fractions $K$ containing the $p$-th roots of unity.
This work is about semi-stable models of $p$-cyclic covers 
of the projective line $\PK$.
It is a continuation of a series of papers on this topic that began with
work of Coleman and McCallum who used rigid geometry in \cite{Co} and \cite{Co-Mc}. 
Our approach is closer to Raynaud's ideas introduced in \cite{Ra1}.
Both \cite{Le1} and \cite{Ma1} dealt with computational aspects of the
semi-stable models of such a cover $C \la \PK$ in the case 
that the branch locus $B$ consists of $K$-points and
has equidistant geometry. The latter means that $\PK$ has a smooth $R$-model
with the property that the points of $B$ have pairwise distinct specializations
on the special fiber. 
In previous work little attention has been given to the determination 
of the finite monodromy, i.e.\ the minimal extension
$K'/K$ such that $C_{K'}$ has a semi-stable model. Let $l \ge 3$ be a prime different from $p$.
It is known (cf.\ \cite{Des} Thm.5.15) that the extension $K'/K$ is 
a certain subextension of
the minimal
extension over which the $l$-torsion points of the jacobian $J(C)$ are rational. 
The tame part of $K'/K$ has been described in a more general setting by Liu and
Lorenzini in \cite{Li-Lo}.
In \cite{Si-Za1} and \cite{Si-Za2} Silverberg and Zarhin study
the finite monodromy groups which can occur for abelian varieties over
a local field of residue characteristic $p\geq 0$. In the case of
abelian surfaces they are able to list these groups in the equal
characteristic $p>0$ case.

So far, in mixed characteristic, there are few techniques available to study the wild finite monodromy, i.e.\
the extension corresponding to the Sylow $p$-subgroup of $\Gal (K'/K)$. 
The literature contains mostly cases where $\deg ([K':K])$ is strictly divisible by $p$.
An exception to this is the work of Krir concerning modular curves (cf.\ \cite{Kr} Thm.2).
This paper is concerned with the wild finite monodromy for $p$-cyclic covers 
in mixed characteristic $(0,p)$.

In section \ref{1}, we review the definitions of the semi-stable models
we shall be concerned with. Associated with $C/K$ are the finite 
monodromy and related Galois groups which will be discussed briefly.

Section \ref{taylor} is the technical heart of the paper. Here, we introduce
the monodromy polynomial ${\cal L}(Y)$.
If the cover is given birationally by the Kummer equation $Z^p=f(X_0)$
then ${\cal L}(Y)$ is related to the logarithmic derivative   
$f'(Y)/f(Y)$. The definition of ${\cal L}(Y)$ is key to the construction
of the stable model as given in Theorem \ref{algo}.
The novelty with respect to the previous work of the authors is that the construction of the 
semi-stable model we give here is very well adapted to study the finite monodromy.
In \cite{Ma1} a polynomial with similar properties as ${\cal L}(Y)$ was defined.
The zeros of either polynomial are centers for the blowups yielding a stable model of $C$.
Yet the polynomial given in \cite{Ma1} has very high degree which allows it to be used for
monodromy purposes only in cases of small genus. 
 
Next, in section \ref{33}, we provide the needed results on the degeneration
of $\mu_p$-torsors. 
The way they are stated here is with a view
towards the later use in conjunction with the monodromy polynomial
and equidistant geometry of the branch locus.
Remark \ref{centers} explains this connection in some detail and we will often
refer to it later. 

Section \ref{44} contains the main result, Theorem \ref{algo}.
It is a new characterization of semi-stable models of 
$p$-cyclic covers of the affine line using the monodromy polynomial
Corollary \ref{df} is an important step towards understanding the finite
monodromy extension, i.e.\ the minimal extension over which $C$ has a stable
model. 

The Galois group of the finite monodromy acts on the special fiber 
$C_{R'} \otimes _{R'} k$ of the stable model:

$$\Gal(K'/K) \hookrightarrow \Aut(C_{R'} \otimes _{R'} k)$$

The bounds derived from \cite{Le-Ma2} Theorem 1.1 for the right hand side 
also limit the possible size of $\Gal(K'/K)$. Hence one can ask if  
$v_p(|\Gal(K'/K)|)$ can be maximal?
(Here $v_p$ is the $p$-adic valuation with $v_p(p)=1$).
In section \ref{66} we give some examples concerning this question for covers that 
are defined
over $\Q_p^{\mbox{\rm \scriptsize tame}}$. The reason for working over this ground
field is that we are mostly interested in the wild part of the monodromy.
We describe various situations where the above can be affirmatively answered.
(cf.\ Examples \ref{gud} and \ref{elli}). The proofs for this section and more on
maximal monodromy will appear in \cite{Le-Ma3}.

\section{Stable models and finite monodromy} \label{1}

Let $R$ be a complete mixed characteristic $(0,p)$
discrete valuation ring with field of fractions $K$ and algebraically closed
residue field $k$. Denote by $v$ the valuation defined by $R$ on $K$ and by $\pi$ a uniformizer.
We will always assume that $R$ contains a primitive $p$-th root of 
unity $\zeta$ and define $\lambda=\zeta-1$. Then $v(\lambda^{p-1})=v(p)$.
Let $C \la \PK$ be a $p$-cyclic cover of smooth projective geometrically irreducible $K$-curves
with genus$(C) \ge 2$.
We assume that the branch locus $B$ of the cover consists of $K$-points,
write $m=|B|-1$, and assume that $B$ has 
equidistant geometry in the sense of the following definition.

\begin{defi} \rm
We say that $B$ has {\it equidistant geometry} if there exists
a smooth $R$-model $\mbox{\rm Proj}(R[u,w])$ for $\PK$ such that the 
points of $B$
have distinct specializations on the special fiber of the model.
\end{defi} 

We can assume that, with respect to the coordinate
$X_0:=u/w$, the cover is given birationally by the equation $Z_0^p=f(X_0)$
with $f(X_0) \in R[X_0]$ monic, $n=\mbox{\rm deg}(f(X_0))$, $(n,p)=1$ and
$n \le m(p-1)$.
Then any two distinct zeros of $f(X_0)$ specialize to distinct elements in $k$.
We denote by $K(C)/K(X_0)$ the function field extension corresponding to
the cover $C \la \PK$. The notion of an equidistant branch locus already
appears in the literature. We refer to \cite{Ra2}, \cite{Gr-Ma} and \cite{Ma1}
for more on this. 

\subsection{Stably marked models and finite monodromy}

For the basic definitions and generalities concerning semi-stable models, we refer to the 
article of Abbes \cite{Ab}.

The following result is due to Deligne and Mumford 
(cf.\ \cite{De-Mu} Cor.2.7 and \cite{Des} Prop.5.7, Lemme 5.16) and is true also in a more 
general setting.

\begin{theorem}[Deligne-Mumford] \label{DM} With the above notation there exists
a minimal extension $K'/K$ such that $C_{K'}$ has a stable model 
$C^{\circ}_{R'}$ over the integral closure $R'$ of $R$ in $K'$. Further $K'/K$
is Galois and acts faithfully on the special fiber:

\begin{equation}
\mbox{\rm Gal}(K'/K) \hookrightarrow \mbox{\rm Aut}_k(C^{\circ}_{R'}\otimes _{R'}k)
\end{equation}
\end{theorem}

\noindent In accordance with Raynaud (cf.\ \cite{Ra2} \S 4) we make the following definition.

\begin{defi} \label{fimo} \rm Let $K'/K$ be the minimal extension over which $C$ has a stable 
model.
We call $\mbox{\rm Gal}(K'/K)$ the {\it finite 
monodromy group} associated to $C$. The extension $K'/K$ is
called the {\it finite monodromy extension}.
$\mbox{\rm Gal}(K'/K)$ has a unique $p$-Sylow subgroup which we call the 
{\it wild monodromy group} $\mbox{\rm Gal}(K'/K)_w$.
The field extension $K'/K'^{\mbox{\rm \scriptsize Gal}(K'/K)_w}$ is called the 
{\it wild monodromy extension}.
\end{defi}

\noindent Keeping the above notation, the following is a consequence of \cite{Liu1} Cor.2.20.
 
\begin{proposition} \label{liu3} Let $K'/K$ be the finite monodromy extension as 
in Definition \ref{fimo}.
Then there exists a semi-stable model $C_{R'}$, defined over the integral closure $R'$ of $R$ in $K'$
minimal with the property that the points in the ramification locus of $C_{K'} \la \PKK$ 
specialize to distinct smooth points on the special fiber.
Further the finite monodromy group operates faithfully on its special fiber:

\begin{equation}\label{inj2}
\mbox{\rm Gal}(K'/K) \hookrightarrow 
\mbox{\rm Aut}_k(C_{R'}\otimes _{R'}k)
\end{equation}

\end{proposition}

\begin{defi} \label{22} \rm The semi-stable model $C_{R'}$, given by Proposition 
\ref{liu3} is called
{\it the stably marked model}.
We refer to the component to which the ramification points specialize as 
{\it original component}.
Any model obtained from the stably marked model by a base change to a ring $R'' 
\supset R'$ we refer to as {\it a stably marked model}.

\end{defi}

The next result is due to Raynaud (cf.\ \cite{Ra1} Thm.2). 
For an effective proof we refer to \cite{Ma1}, Thm.3.2.2.

\begin{proposition} \label{tree}
Let $C \la \PK$ be as above. Then 
the dual graph of the special fiber of a stably marked model  
of $C$ is an oriented tree whose origin corresponds to the original component (cf.\ Def.\ref{22}).
(For short we say that the special fiber is \it{tree-like}).
\end{proposition}

\begin{proposition} \label{liu}
Let $C \la \PK$ be as above. Assume that for each irreducible component $E$ of
genus $>0$ in the stable reduction there exists a flat quasi-projective $R$-model 
${\mathcal C}$ for
$C$, such that ${\mathcal C} \otimes _R k$ is birational to $E$.
Then $C$ has a stable model over $R$.
\end{proposition}

\begin{proof}
By Proposition \ref{tree} we know that the special fiber is tree-like.
Now the claim is an immediate consequence of \cite{Liu2} chap.10, Prop.3.44.
\end{proof}

As stated before, we shall always assume that the branch locus of $C \la \PK$ 
has equidistant geometry and consists of rational points. 
This has consequences for the image of the injection
\eqref{inj2}. Namely, any element of its image will have trivial action on the
original component (cf.\ Definition \ref{22}) 
of the stably marked model. The group of such automorphisms
we denote by $\mbox{\rm Aut}_k(C_{R'}\otimes _{R'}k)^{\circo}$. Hence we have

\begin{equation}\label{inj}
\mbox{\rm Gal}(K'/K) \hookrightarrow 
\mbox{\rm Aut}_k(C_{R'}\otimes _{R'}k)^{\circo}
\end{equation}

\begin{defi} \label{mamo} \rm If the injection \eqref{inj} is surjective
we say that $C$ has {\it maximal monodromy}.
If $v_p(|\Gal (K'/K)|) = v_p(|\mbox{\rm Aut}_k(C_{R'}\otimes _{R'}k)^{\circo}|)$
we say that $C$ has {\it maximal wild monodromy}.
\end{defi}

The following section is the technical heart of the paper and introduces
the monodromy polynomial.

\section{The monodromy polynomial} \label{taylor}
For this section let $R$ be a discrete valuation ring as defined in section \ref{1} and 
$f(X_0) \in R[X_0]$ a
monic polynomial of degree $n$ prime to the residue characteristic $p$
of $R$. Let $m$ be the number of distinct zeros of $f(X_0)$ in an algebraic closure of
$\mbox{\rm frac}(R)$ and $r$ the greatest integer such that $rp<n$.
For $X_0=X+Y$ Taylor expansion yields
\begin{equation} \label{tt}
f(X_0)=f(X+Y)=s_0(Y)+s_1(Y)X+s_2(Y)X^2+ \dots +s_n(Y)X^n \in R[X,Y]
\end{equation}
which we view as a polynomial in two variables. 

\begin{defi}\rm \label{deco}
If $f(X_0)$ as above can be written
\begin{equation} \label{exp}
f(X+Y)=s_0(Y)\left(H(X,Y)^p-\sum_{i=r+1}^{n}A_i(Y)X^i\right)
\end{equation}
with $H(X,Y)=1+a_1(Y)X+a_2(Y)X^2+ \dots +a_r(Y)X^r \in K(Y)[X]$, $A_i(Y) \in K(Y)$
and $r$ as introduced earlier,
we call this a {\it special decomposition} of $f(X_0)$.
\end{defi} 

The existence and properties of such a decomposition are established in
Lemma \ref{H}.
We set 
$S_1(Y)=s_1(Y)/\mbox{\rm gcd}(s_0(Y),s_1(Y))$ and
$S_0=s_0(Y)/\mbox{\rm gcd}(s_0(Y),s_1(Y))$.
Then $S_1(Y), S_0(Y) \in R[Y]$, $f'(Y)/f(Y)=S_1(Y)/S_0(Y)$, $\deg (S_1(Y))=m-1$, $\deg (S_0(Y))=m$  
and $(S_0(Y),S_1(Y))=1$.
Also, for $t \in \mathbb N$, we will use the following notation $${\frac{1}{p} \choose t}
=\frac{\prod_{i=0}^{t-1}\left(\frac{1}{p}-i\right)}{t!}.$$
For later use we note that the $p$-adic value of these binomial coefficients is given by

\begin{equation} \label{value}
-v_p({\frac{1}{p} \choose t})=t+[t/p]+[t/p^2]+ \dots 
\end{equation}
\noindent where $v_p(p)=1$.

\begin{lemma} \label{int}
Let $n$ be a positive integer not divisible by $p$, and $r$
the greatest integer, such that $rp<n$. Consider the equivalence relation
$\mathcal R$ on $M=\{1,2, \dots ,n\}$ that identifies $a$ with $b$ if 
there exists an integer $e \ge 0$ such that $a=p^eb$ or $b=p^ea$.

Then the set $M_0=\{r+1,r+2, \dots, n\}$ is a minimal complete system of
representatives for the quotient $M/\mathcal R$. In particular $M_0$
contains exactly one power $p^{\alpha}$ of $p$, characterized by $p^{\alpha} < n <
p^{\alpha +1}$.

\end{lemma}

\noindent The proof is elementary and omitted. 
 
\begin{lemma} \label{H} i) A special decomposition as given by equation \eqref{exp}
above exists and is unique.

\noindent ii) There exist $t_i(Y)$, $T_i(Y) \in R[Y]$ such that
$$a_i(Y)={\frac{1}{p} \choose i}\frac{s_1(Y)^i+pt_i(Y)}{s_0(Y)^i} 
={\frac{1}{p} \choose i}\frac{S_1(Y)^i+pT_i(Y)}{S_0(Y)^i} 
\quad \mbox{for} \quad  1 \le i \le r $$ 
where $t_1(Y)=0$ and for $1<i\le r$, $\deg (t_i(Y)) = i(n-1)$. 

\noindent iii) There exist $c_i \in K$, $T(Y) \in R[Y]$ and $N_i(Y) \in R[Y]$ monic such that
$$A_i(Y)=c_i\frac{N_i(Y)}{s_0(Y)^i} \quad \mbox{ for } r+1 \le i \le n \quad
\mbox{and} $$

$$A_{p^\alpha}(Y)=-{\frac{1}{p} \choose p^{\alpha-1}}^p \frac{S_1(Y)^{p^\alpha}+pT(Y)}
{S_0(Y)^{p^{\alpha}}} \quad
\mbox{with} \quad v_p(c_{p^\alpha})=v_p({\frac{1}{p} \choose p^{\alpha-1}}^p) \le 0$$

\noindent where $p^\alpha$ is the power determined by Lemma \ref{int} and for $\alpha=0$ we set
${\frac{1}{p} \choose \frac{1}{p}}=1$.
 
\end{lemma}

\begin{proof}
i) First observe that for $r=0$ we can take $H(X,Y)=1$ and then the statement is 
true. The case of $r>0$ remains to be considered.
One computes the functions $a_i(Y)$ recursively from equation \eqref{exp} 
starting
with $a_1(Y)=s_1(Y)/(ps_0(Y))$. This requires only solving linear equations and
therefore the existence and uniqueness are obvious.

For later use we shall show the following property: $S_0(Y)^ia_i(Y) \in K[Y]$ for
$1 \le i \le r$. We proceed by induction on $i$. For $i=1$ this is immediate from the
expression for $a_1(Y)$ given above. We therefore assume that the property holds
for all $i$ with $1 \le i \le l<r$. From \eqref{exp} we get
$$f(X+Y)=s_0(Y)(1+a_1(Y)X+ \dots + a_l(Y)X^l)^p+ps_0(Y)a_{l+1}X^{l+1} \quad (X^{l+2}).$$
Comparing the coefficients of $X^{l+1}$ on each side we get an expression for $a_{l+1}(Y)$. 
Now taking into account that $S_0(Y)^{l+1}s_{l+1}(Y)/s_{0}(Y) \in K[Y]$ the induction
hypothesis yields $S_0(Y)^{l+1}a_{l+1}(Y) \in K[Y]$.

ii) To derive these properties, an alternate approach, such as that utilized under i) is more advantageous.
From \eqref{exp} we get the following identity in $K(Y)[[X]]$:
$$H(X,Y)=\left(\frac{f(X+Y)}{s_0(Y)}-\sum_{i=r+1}^{n}A_i(Y)X^i\right)^{1/p}$$
$$=\left(1+\frac{s_1(Y)}{s_0(Y)}X+\frac{s_2(Y)}{s_0(Y)}X^2 +\dots +\frac{s_n(Y)}{s_0(Y)}X^n
-\sum_{i=r+1}^{n}A_i(Y)X^i\right)^{1/p}
=(1+F(X))^{1/p}$$ where
$$F(X)=\sum_{i=1}^{n}\frac{s_i(Y)}{s_0(Y)}X^i-\sum_{i=r+1}^{n}A_i(Y)X^i \in XK(Y)[[X]].$$
We can develop this in $K(Y)[[X]]$ to get $(1+F(X))^{1/p}$ with value $1$ at $X=0$:
$$1+{\frac{1}{p} \choose 1}F(X)+{\frac{1}{p} \choose 2}F(X)^2 + 
{\frac{1}{p} \choose 3}F(X)^3 + \dots$$
Because $\mbox{deg}_X(H(X,Y)) \le r$, we can compute modulo $X^{r+1}$ and
therefore only need to evaluate the following expansion

$$1+{\frac{1}{p} \choose 1}F_0(X)+{\frac{1}{p} \choose 2}F_0(X)^2 + 
{\frac{1}{p} \choose 3}F_0(X)^3 + \dots
\quad \mbox{ where } F_0(X)=\sum_{i=1}^{n}\frac{s_i(Y)}{s_0(Y)}X^i.$$
From this one reads the coefficients of $H(X,Y)$. Notice that by \eqref{value}
$${\frac{1}{p} \choose t}^{-1}{\frac{1}{p} \choose s} \in pR
\quad \mbox{for any two integers } t>s.$$ 
This justifies the binomial coefficient in the first expression for $a_i(Y)$.
Under i) we have seen that $S_0(Y)^ia_i(Y) \in K[Y]$. Further, $s_1(Y)/s_0(Y)=S_1(Y)/S_0(Y)$
which allows us to write the second expression given for $a_i(Y)$
with $T_i(Y)$ a polynomial. 

iii) Now that an expression for $H(X,Y)$ is determined we can use \eqref{exp}
to compute the $A_i(Y)$. The expression $A_i(Y)=c_i \frac{N_i(Y)}{s_0(Y)^i}$ 
is then immediate. To obtain $A_{p^{\alpha}}(Y)$
observe that $p^{\alpha-1} < r+1 < p^\alpha$. Therefore

\begin{equation} \label{new}
A_{p^\alpha}(Y)=\frac{s_{p^\alpha}(Y)}{s_0(Y)}-a_{p^{\alpha-1}}(Y)^p
+\mbox{other terms}
\end{equation}

Using \eqref{value} it is straightforward to verify that the {\it other terms}
and $s_{p^\alpha}(Y)$
all have higher divisibility by $p$ than $a_{p^{\alpha-1}}(Y)^p$.
Like for the $a_i(Y)$ this allows us to write

$$A_{p^\alpha}(Y) = - {\frac{1}{p} \choose p^{\alpha-1}} 
\frac{s_1(Y)^{p^\alpha} + p t_{p^\alpha}(Y)}{s_0(Y)^{p^\alpha}} \quad
\mbox{\rm with } t_{p^\alpha}(Y) \in R[Y]$$

We proceed to get the expression for $A_{p^\alpha}(Y)$ given in the lemma. 
One can rewrite \eqref{new} as 

$$s_{p^\alpha}(Y) = s_0(Y)[H_{p^\alpha}(Y)+A_{p^\alpha}(Y)]$$ 

\noindent where $H_{p^\alpha}(Y)$ is the coefficient of $X^{p^\alpha}$
in $H(X,Y)^p$.
Notice that 
$S_0(Y)^{p^\alpha}H_{p^\alpha}(Y) \in K[Y]$ by what was said under i).
We hence conclude that $S_0(Y)^{p^\alpha}A_{p^\alpha}(Y) \in K[Y]$.
This yields the term $T(Y)$ in the formula for $A_{p^\alpha}(Y)$.
\end{proof}

\begin{defi} \label{lpd} \rm In the situation of Lemma \ref{H} part iii) we
define ${\mathcal L}(Y)=S_1(Y)^{p^\alpha}+pT(Y)$.
It is a polynomial of degree $p^{\alpha}(m-1)$ and we call it the 
{\it monodromy polynomial} of $f(Y)$.
\end{defi}

\noindent From the definition one immediately obtains 
${\mathcal L}(Y) \equiv S_1(Y)^{p^\alpha}$ mod $p$. Notice that $f'(Y)/f(Y)=S_1(Y)/S_0(Y)$.

\section{Degeneration of $\mu_p$-torsors} \label{33}

In this section we have collected various results on the degeneration 
of $\mu_p$-torsors. These are not new and similar results have appeared already 
in various places.
For this reason we do not include proofs here but refer to the existing literature.

\subsection{The general case} \label{here}
Let $R$ be a discrete valuation ring as in section \ref{1} and recall the definition
of $\llll$ given there. Let $C \la \PK$ be a $p$-cyclic cover of smooth projective 
geometrically irreducible $K$-curves.
Unlike in the rest of the paper, in section \ref{here} we make no 
assumptions about the branch locus of the cover $C \la \PK$.

\begin{defi} \label{gauss} \rm
The valuation on $K(X_0)$ corresponding to the DVR $R[X_0]_{(\pi)}$ is called
the {\it Gauss valuation} $v_{X_0}$ with respect to $X_0$.
\end{defi}

$$ \mbox{\rm We then have} \quad 
v_{X_0}\left( \sum_{i=0}^m a_iX_0^i \right) 
= \mbox{\rm min}\{v(a_i) | 0 \le i \le m\}.$$
Notice that a change of coordinates $X=(X_0-y)/\rho$ for $\rho, y \in R$ gives
rise to a Gauss valuation $v_X$. It is equivalent to $v_{X_0}$
iff $\rho$ is a unit.
Further, these
valuations are exactly those that come from the local rings at generic points
of components in the semi-stable models for $\PK$.

\begin{lemma} \label{aprox}
Let $C \la \PK$ be a $p$-cyclic cover of curves and
$\mbox{\rm Proj}(R[u,w])$, with $X=u/w$, a smooth $R$-model for $\PK$.
Denote by ${\mathcal C}_R$ its normalization in $K(C)$.
Suppose the cover is given birationally by an equation $Z^p=f(X)$ with
$f(X) \in R[X]$ and such that $\pi$ does not divide $f(X)$.
Then the cover of special fibers ${\mathcal C}_R\otimes_Rk 
\la \mbox{\rm Proj}(k[u,w])$ is separable and irreducible iff
$$\mbox{\rm max}\{v_X(h(X)^p-f(X)) \mbox{ for } h(X)\in R[X]\}=v(\lambda^p).$$
\end{lemma}

\noindent The proof is elementary and essentially that of \cite{Le1} Prop.1 
or \cite{Gr-Ma} III.Prop.1.1.

\subsection{The equidistant case}
In this section we work under the hypotheses of section \ref{1}, i.e.\
let $C \la \PK$ be a $p$-cyclic cover and assume that the  
branch locus $B$ consists of $K$-points having equidistant geometry
with respect to the coordinate $X_0$.
We suppose that the cover is given by $Z_0^p=f(X_0)$ where $X_0=u/w$ and
with $f(X_0) \in R[X_0]$ monic of degree $n$ prime to $p$. 
We will also use the notation introduced in Definition \ref{deco}. 

\begin{figure}[h]
\begin{center}
\input{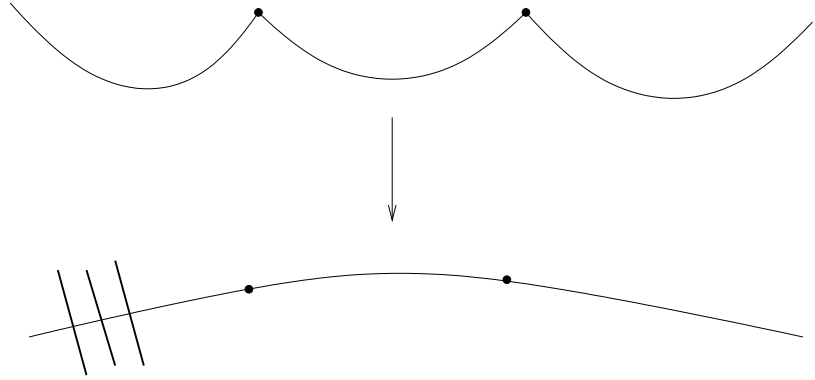}
\caption{${\mathcal C}_{R'} \otimes _{R'} k \la \Pk$ with singularities and branch locus}
\label{f1}
\end{center}
\end{figure}

\begin{remark} \label{centers} \rm
a) In the above situation, consider the 
flat projective $R$-model ${\mathcal C}_R$ for $C$ 
obtained by normalizing $\mbox{\rm Proj}(R[u,w])$ in the function field $K(C)$.
From ${\mathcal C}_R$ we obtain the stably marked model $C_{R'}$ for $C$ by a series
of blowups, after passing to a finite extension $R'$. 
The uniqueness of $C_{R'}$ implies that the $p$-cyclic group action on $C$ extends to 
$C_{R'}$ and the quotient by this action is a
semi-stable model for ${\mathbb P}^1_{K'}$, where $K'=\mbox{\rm frac}(R')$ 
(cf.\ \cite{Ra1} appendice).
This semi-stable model is obtained from $\mbox{\rm Proj}(R'[u,w])$ by a series of blowups.
By construction, these are all centered on the affine patch $\mbox{\rm Spec}(R'[X_0])$
and hence correspond to ideals of the form $(X_0-y,\rho) \subset R'[X_0]$. 
We call $y$ the center and
$\rho$ the radius, and often speak of the Gauss valuation $v_X$ with $X=(X_0-y)/\rho$
instead.
All blowups we will consider are of this type and on 
$\mbox{\rm Proj}(R^{\scriptsize alg}[u,w])$.
Blowing up $(X_0-y,\rho) \subset R'[X_0]$
gives an exceptional divisor which in turn, 
by normalization in $K'(C)$, yields irreducible components of the stably marked model.
The DVRs at their generic points are the ones that
dominate the DVR of $v_X$. An interesting case, when there is only one, is
treated in Proposition \ref{deg}. 

A local computation shows that the centers specialize to points below
the singularities (cusps) of ${\cal C}_R \otimes _R k$ and are outside the 
specialization of the branch locus.
Their $X_0$-coordinates are zeros of $\bar{f}'(X_0)$, hence also zeros of
$\bar{S}_1(X_0)=\bar{f}'(X_0)/\mbox{\rm gcd}(\bar{f}'(X_0),\bar{f}(X_0))$.
In particular, if $X_0=y$ is such a center, then $y\in R^{\scriptsize alg}$ and
$\bar{f}(\bar{y}) \in  k-\{0\}$. 
A schematic picture of this situation is given in figure \ref{f1}.

b) It is immediately derived from the equation that the cover of special fibers
${\mathcal C}_R \otimes _R k \la \Pk$ is purely inseparable. Therefore,
to yield a component of genus $>0$ in the stable reduction of $C$, a blowup
(on $\mbox{\rm Proj}(R[u,w])$)
must have a radius $\rho$ with $v(\rho)>0$.

\end{remark}

\noindent With the above notation we obtain the following.

\begin{proposition}\label{deg}
The Gauss valuation corresponding to $X=(X_0-y)/\rho$ for $\rho, 
y \in R^{\scriptsize alg}$,
$\bar{f}(\bar{y})\not= 0$
induces a separable and irreducible component in a semi-stable model of $C$ iff
the following two conditions hold:

\noindent i) $$H(\rho X,y) \in R^{\scriptsize alg}[X]$$
 
\noindent ii) $$v_X\left( \frac{f(\rho X+y)}{s_0(y)}-H(\rho X,y)^p\right)=v(\lambda^p).$$\\

\noindent In this case, and with the change of coordinates 
$Z_0/(s_0(y)^{1/p})=\lambda Z_1+H(\rho X,y)$, the equation $Z_0^p=f(X_0)$ yields

$$
Z_1^p+\frac{p}{\lambda}Z_1^{p-1}H(\rho X,y)+ 
\frac{{p \choose 2}}{\lambda^2}Z_1^{p-2}H(\rho X,y)^2+ \dots $$
\begin{equation} \label{degeq}
+\frac{p}{\lambda^{p-1}}Z_1H(\rho X,y)^{p-1}
-\frac{1}{\lambda^p}
\sum_{i=r+1}^{n}A_i(y)\rho^iX^i=0
\end{equation}

\noindent which, for $v(\rho)>0$, 
in reduction gives the irreducible Artin-Schreier equation
\begin{equation} \label{arti}
Z_1^p-Z_1-\left(\frac{1}{\lambda^p}\sum_{i=r+1}^{n}
A_i(y)\rho^iX^i\right)^-=0
\end{equation}
\noindent defining an \'etale cover of ${\mathbb A}^1_k$.
\end{proposition}

\begin{proof}
This is a straightforward translation of Lemma \ref{aprox} into the terminology
of Definition \ref{deco}. We therefore omit the details of the proof.
\end{proof}



\section{Construction of a stably marked model} \label{44}

Now we can state the theorem describing a stably marked model. 
We work under the hypotheses of section \ref{1}.
Suppose that the cover is given birationally by $Z_0^p=f(X_0)$ with
$f(X_0)$ monic of degree $n$ prime to $p$ and such that, with respect to $X_0$,
the branch points have distinct specializations. 
Let $r$ be the greatest
integer such that $rp<n$ and let $\alpha$ be determined by $p^{\alpha} <n< p^{\alpha +1}$.
We write $|B|-1=m$ and can assume that
$f(X_0)$ has been chosen such that $m(p-1) \ge n$. Choosing $f(X_0)$ this way limits the size 
of ${\mathcal L}(Y)$. We point out that, due to the freedom in choosing $f(X_0)$, the polynomial
${\cal L}(Y)$ is not uniquely determined by the cover $C \la \PK$.

\begin{theorem} \label{algo}
The components of genus $>0$ of the stably marked model of $C$
correspond bijectively to the Gauss valuations $v_{X_j}$
with $\rho_j X_j=X_0-y_j$ where $y_j$ is a zero of the monodromy polynomial
 ${\mathcal L}(Y)$ and $\rho_j \in R^{\mbox{\rm \scriptsize alg}}$ is
such that
$$ v(\rho_j)=\mbox{\rm max}\{ 
\frac{1}{i}v\left(\frac{\lambda^p}{A_i(y_j)}\right) 
\mbox{ for } r+1 \le i \le n \}.$$

\noindent The dual graph of the special fiber of the stably marked model of $C$ is an oriented tree
whose ends are in bijection with the components of genus $>0$.
\end{theorem}

\begin{proof} 
The fact that the special fiber is tree-like is a consequence of Proposition
\ref{tree}. The statement about the ends follows from \cite{Ra2} 2.4.8.
Let $y_j$ be a zero of ${\mathcal L}(Y)$. 
Then $\overline{S_1(y_j)}=0$ and hence $y_j \in 
R^{\rm \scriptsize alg}$. As $S_1(Y)$ and $s_0(Y)$ are relatively prime we conclude that
$s_0(y_j)=f(y_j)$ is a unit (cf. Remark \ref{centers}).
We will show that for $\rho_j$
as defined above, one obtains a component of genus $>0$.
As a first step, we need to show that $v(\rho_j) > 0$. Therefore assume $v(\rho_j) \le 0$.
Then we get $v(A_i(y_j)) \ge v(\lambda^p)$ for all $r+1 \le i \le n$. Using equation
\eqref{exp} of Definition \ref{deco}, we conclude that $H(X,y_j) \in R[y_j][X]$ and therefore
equation \eqref{arti} in Proposition \ref{deg} applies to show that the Gauss 
valuation corresponding to
$\rho=1$ and $y=y_j$ induces a separable (not necessarily irreducible)
component in a semi-stable reduction of $C$.
Now this Gauss valuation is $\mbox{PGL}_2(R[y_j])$-equivalent to the one given by $X_0$,
so both induce the same component in the semi-stable reduction. This is a contradiction,
as we know that the latter induces a ${\bf \mu}_p$-torsor. We conclude $v(\rho_j) > 0$.
Next
consider the equation \eqref{degeq} in Proposition \ref{deg} for $y=y_j$ and $\rho=\rho_j$:
$$
Z_1^p+\frac{p}{\lambda}Z_1^{p-1}H(\rho_j X,y_j)+ 
\frac{p}{\lambda^2}Z_1^{p-2}H(\rho_j X,y_j)^2+ \dots $$
$$-Z_1H(\rho_j X,y_j)^{p-1}-\frac{1}{\lambda^p}\sum_{i=r+1}^{n}A_i(y_j)\rho_j^iX^i=0.
$$   

Then, by definition of $\rho_j$, the polynomial 
$\frac{1}{\lambda^p}\sum_{i=r+1}^{n}A_i(y_j)\rho_j^iX^i$ is in $R[y_j,\rho_j][X]$
and $s_0(y_j)$ is a unit as the branch locus has equidistant geometry.
Now the relation \eqref{exp} implies the same for $H(\rho_j X,y_j)$.
Altogether we see that the above equation is in $R[y_j,\rho_j][X,Z_1]$.
Now, taking into account Proposition \ref{deg}, we conclude that in reduction
we get an irreducible Artin-Schreier equation.
Further, using Artin-Schreier theory, Lemma \ref{int} 
implies that the equation 
cannot have genus zero, because by definition of ${\mathcal L}(Y)$ we have
$A_{p^\alpha}(y_j)=0$. This fact is the motivation for the definition of
${\mathcal L}(Y)$.

Conversely, suppose the Gauss valuation corresponding to $X=(X_0-y)/\rho$ 
induces
a component of genus $>0$ in the stable reduction of $C$.
Then, by Remark \ref{centers}, $y \in R^{\mbox{\rm \scriptsize alg}}$, $v(\rho) > 0$ and 
$s_0(y) \in R^{\mbox{\rm \scriptsize alg}}$ is a unit.
By Proposition \ref{deg} an Artin-Schreier equation for
the component will be given by reducing
equation \eqref{degeq} to $k$. 
By assumption the conductor of this equation is strictly bigger than $1$
and in turn corresponds by Lemma \ref{int} to a unique integer $t$ with
$r+1 \le t \le n$ and $t\not=p^\alpha$.
Now for any ${y}_1$ with 
$v({y}_1-y)\ge v(\rho)$, the
Gauss valuation corresponding to ${X}_1=(X_0-{y}_1)/\rho$ yields
the same component in the stable reduction. In particular the associated
Artin-Schreier equation has to have the same conductor, and hence the same
value of $t$, as defined above.
We conclude that $v(A_t(y))=v(A_t({y}_1))$ and $v(A_t(y)\rho^t)=v(\lambda^p)$
for any ${y}_1$ 
with $v({y}_1-y)\ge v(\rho)$.
Using Lemma \ref{H} iii) and observing that $s_0(y)$ and $s_0({y}_1)$ are units
we get  $v(N_t(y))=v(N_t({y}_1))$.

This implies that for all zeros $z$ of $N_t(Y)$ we have $v(z-y) < v(\rho)$.
In particular there exist a $\rho' \in R^{\mbox{\rm \scriptsize alg}}$ such that
$v(z-y) < v(\rho') < v(\rho)$ as $N_t(Y)$ has only finitely many zeros.
Furthermore we can assume that $\rho'$ has been chosen such that
for all $\tilde{y}$ with $v(\rho) > v(y-\tilde{y}) > v(\rho')$
the point $\tilde{y}$ is not a center for a component of genus $>0$ in the 
stable reduction of $C$. We fix $\tilde{y}$ such that $v(\rho) > v(y-\tilde{y}) > v(\rho')$.
Now pick $\tilde{\rho} \in R^{\mbox{\rm \scriptsize alg}}$ such that
$$v(\tilde{\rho})=\mbox{\rm max}\{ 
\frac{1}{i}v\left(\frac{\lambda^p}{A_i(\tilde{y})}\right) 
\mbox{ for } r+1 \le i \le n \}.$$
By choice of $\rho'$ we have $v(z-\tilde{y}) < v(\rho')$ for all zeros $z$ of $N_t(Y)$.
Therefore, as $v(y-\tilde{y}) > v(\rho')$, we get $v(N_t(y))=v(N_t(\tilde{y}))$.
Now the definition of $\tilde{\rho}$ yields $v(\tilde{\rho}) \ge v(\rho)$.
Also by construction equation \eqref{degeq} of Proposition \ref{deg} 
for $y=\tilde{y}$ and $\rho=\tilde{\rho}$
will have integral coefficients and be irreducible in reduction.
We already know that the corresponding component in a semi-stable reduction
must have genus zero.
Therefore $v(A_t(y))=v(A_t(\tilde{y}))$ even implies $v(\tilde{\rho}) > v(\rho)$.
Further, using Lemma \ref{H} iii), genus zero implies 
$$v({\frac{1}{p} \choose p^{\alpha-1}}^p{\mathcal L}(\tilde{y})\tilde{\rho}^{p^\alpha})=v(\lambda^p).$$
By assumption the component corresponding to $X=(X_0-y)/\rho$ has genus $>0$
consequently
$$v({\frac{1}{p} \choose p^{\alpha-1}}^p{\mathcal L}(y) \rho^{p^\alpha}) \ge v(\lambda^p).$$
Altogether one concludes that 
$$v({\mathcal L}(\tilde{y})) < v({\mathcal L}(y)) \quad \mbox{\rm for all } \tilde{y}
\mbox{ with} \quad v(\rho') < v(y-\tilde{y}) < v(\rho).$$
We claim this implies that there exists a zero $z$ of ${\mathcal L}(Y)$ with
$v(z-y) \ge v(\rho)$. To see this assume the contrary, namely $v(z-y) < v(\rho)$ for all zeros $z$ of
${\mathcal L}(Y)$. Chose $\tilde{y}$ such that $v(\rho') < v(y-\tilde{y}) < v(\rho)$ and 
$v(y-\tilde{y}) > v(z-y)$ for all zeros $z$ of ${\mathcal L}(Y)$. This is possible because
${\mathcal L}(Y)$ has only finitely many zeros. Then $v(z-\tilde{y})=v(z-y+y-\tilde{y})=v(z-y)$ for 
all zeros $z$ of ${\mathcal L}(Y)$ and hence $v({\mathcal L}(y))=v({\mathcal L}(\tilde{y}))$. This is a
contradiction thus proving the claim.

\end{proof}

\begin{figure}[h]
\begin{center}
\input{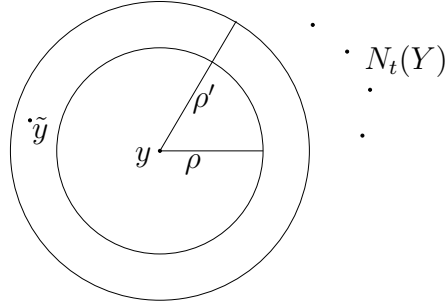}
\caption{Illustrating the $p$-adic geometry in the second part of the 
proof of Theorem \ref{algo}}
\label{3fibers}
\end{center}
\end{figure}

\begin{corollary} \label{df}
Let $K'/K$ be the minimal extension such that $C_{K'}$ has a stable model   
over the integral closure $R'$ of $R$ in $K'$
and denote by F the splitting field of ${\mathcal L}(Y)$
over $K$. Let $E$ be the field obtained from $F$ by adjoining
$\rho_i$ as well $(s_0(y_i))^{1/p}$ for each $i$. 
Then $K' \subseteq E$. In particular $C$ has a stable model over the integral closure 
of $R$ in $E$.
\end{corollary}

\begin{proof} Observe that all the blowups needed to obtain a stably
marked model are defined over $E$, and the sum of the genera of the
exceptional divisors is equal to the genus of the generic fiber
(cf.\ Proposition \ref{tree}). Now, applying Proposition \ref{liu},
the result follows.
\end{proof}


\begin{remark} \label{gg} \rm 
Assume that $C$ has potentially good reduction.
Then the special fiber of the stably marked model $C_{R'} \otimes_{R'} k$
has 2 irreducible components:
the original component and a component $E_k$ with the same geometric genus as $C$.
We denote their unique point of intersection by $\infty$.
By Proposition \ref{deg} $E_k$ is a $p$-cyclic cover of
${\mathbb P}_k^1$ which is ramified only at $\infty$. 
Let $\mbox{\rm Aut}(E_k)_{\infty}$ be the group of automorphisms of $E_k$ 
that leave the point $\infty$ fixed.
Then this group is isomorphic to the group $\mbox{\rm Aut}_k(C_{R'}\otimes _{R'}k)^{\circo}$
which was defined above Definition \ref{mamo}.
Let $q$ be the highest power of $p$ dividing $[K':K]$. From the injection \eqref{inj}
we conclude that $q$ also divides the order of the $p$-Sylow subgroup
of $\mbox{\rm Aut}(E_k)_{\infty}$. Let $m=|B|-1$ and write $m=lp^s+d$ with 
$(l,p)=1$ and $1 \le d \le p-1$. Then using \cite{Le-Ma2} Theorem 1.1
we get the following
\begin{itemize}

\item[1)] If $d\not= 1$ then $|\mbox{\rm Aut}(E_k)_{\infty}|$ is at most divisible by $p$
hence $q \in\{1,p\}$ (cf.\ \cite{Le-Ma2}, Thm.1 II.b).
This explains why in \cite{Le1} Theorem 1 we got
a criterion for good reduction with $[K':K]$ at most divisible by $p$.

\item[2)] If $d=1$ and $l>1$ then $|\mbox{\rm Aut}(E_k)_{\infty}|$ is at most 
divisible by $p^{s+1}$ ($p^{s}$ if $p=2$).
Further the $p$-Sylow subgroup of $\mbox{\rm Aut}(E_k)_{\infty}$ injects into
an extraspecial $p$-group
(cf.\ \cite{Le-Ma2}, Thm.1 II.d and \cite{Su} chap.\ 4 for the notion of extraspecial group).
In this case, by Definition \ref{lpd}, the polynomial ${\mathcal L}(Y)$ has degree $lp^{s+\alpha}$ 
(recall that $\alpha$ is defined by $p^\alpha < n < p^{\alpha+1}$). 

\item[3)] If $d=l=1$ then $|\mbox{\rm Aut}(E_k)_{\infty}|$ is at most 
divisible by $p^{2s+1}$. Again the $p$-Sylow subgroup of 
$\mbox{\rm Aut}(E_k)_{\infty}$ injects into an extraspecial $p$-group
(cf.\ \cite{Le-Ma2}, Thm.1 II.a).
Recall that $n=\deg (f(X_0))$ and $m \le n \le (p-1)m$. 
Using $m=1+p^s$ we get 
$p^s < n < p^{s+1}$ and Definition \ref{lpd} shows that
${\mathcal L}(Y)$ has degree $p^{2s}$.

\end{itemize}

Fix $n$ and $m=p^s+1$.
Point 3), together with Example \ref{gud} below,
show that $\mathcal L(Y)$ has the smallest degree one could expect for 
this $n$ and $m$, as it
has been defined independently of the type of degeneration of the model.

\end{remark}

\begin{remark} \label{ggg} \rm
As explained earlier, there are bounds for the finite monodromy group
$\mbox{Gal}(K'/K)$ derived from the injection \eqref{inj} and the bounds 
for automorphism groups in positive characteristic (cf.\ \cite{Le-Ma2} Thm.1).
Examples \ref{gud} and \ref{elli}
give reason to believe that these bounds on $\mbox{Gal}(K'/K)$  are sharp, as there we have 
maximal wild monodromy. 
Many more examples of this type will be given in forthcoming work of the 
authors (cf. \cite{Le-Ma3}).
Due to the relationship between ${\cal L}(Y)$ and $K'/K$ given by Thm.\ref{algo} and 
Cor.\ref{df}, the degree of ${\cal L}(Y)$ should reflect these bounds, 
if the degree of ${\cal L}(Y)$ is optimal. We fix $n$ and $m$.

\noindent a) Remark \ref{gg} above shows that for $m=1+p^s$ this is the case.

\noindent b) For $m=lp^s+d$, $(l,p)=1$, $1 \le d \le p-1$ : Assume there exists a cover 
$C \la \PK$ whose stable reduction
contains a component $E$ of genus $p^a (p-1)/2$ with $a=[\log_p(m-1)]$ 
maximal. $E$ is a $p$-cyclic cover of $\Pk$ that is ramified in exactly one point that we denote $\infty$. The assumption on the genus of $E$ is hence equivalent to the conductor
of this cover at $\infty$ being equal to $p^a+1$.
Further suppose that $C$ has maximal wild monodromy (cf.\ Def.\ref{mamo}). 
A small calculation shows that, 
with few exceptions that we will not look at here, $a=\alpha$, where $p^{\alpha} < n < p^{\alpha +1}$.
\cite{Le-Ma2} Thm.1 yields  $v_p(|\Aut (E)|) \le 2\alpha +1$ where equality is possible.
On the other hand $\deg ({\cal L}(Y)) = (m-1)p^{\alpha}$ and $p^{\alpha} \le m-1 < p^{\alpha +1}$
imply $p^{2 \alpha} \le \deg ({\cal L}(Y)) < p^{2 \alpha +1}$. 
Comparing with the size of $v_p(|\Aut (E)|)$ we see that $v_p(\deg({\cal L}(Y)))$ is optimal if
one makes the above assumptions on the existence of certain covers with prescribed reduction. 

Notice that the discrepancy of $1$ between the right hand side
of $v_p(|\Aut (E)|) \le 2\alpha +1$ and the left hand side of 
$2 \alpha < v_p(\deg ({\cal L}(Y)))$ is explained by the fact that in
Cor.\ref{df} certain $p$-th roots have to be adjoined.
\end{remark}


We finish by giving some indication on our examples of maximal finite
monodromy over $\Qpt$. A more extensive version of this will be published
elsewhere (cf.\ \cite{Le-Ma3}).

\begin{exa}[Potentially good reduction with $m=1+p^s$] \label{gud} \rm

Let $p>2$, $q=p^n$, $n\geq 1$,
$K={\mathbb Q}_p^\ex(p^{p/(q+1)})$ 
and $C \la \PK$ be given birationally by the equation 

$$Z_0^p=f(X_0)=1+p^{p/(q+1)}X_0^q+X_0^{q+1}.$$

\noindent
Then $C$ has potentially good reduction and
${\cal L}(Y)$ is irreducible over $K$.
The finite monodromy
extension $K'/K$ is obtained from the splitting field of ${\cal L}(Y)$ by adjoining
$f(y)^{1/p}$ for $y$ any zero of ${\cal L}(Y)$.
It is Galois with the extraspecial group
of exponent $p$ and order $pq^2$ (see \cite{Su} ch.\ 4 for these groups). 
Further the finite monodromy is maximal in the sense of Definition \ref{mamo}.
\end{exa}

\begin{exa}[Genus $2$-curves and wild monodromy] \label{elli} \rm

Now we consider the case of $p$-cyclic covers where $p=2$ and 
$m=5$, i.e.\ genus 2 curves over a $2$-adic field in ${\mathbb Q}_2^{{\rm \scriptsize tame}}$.
In this case, there are three possible types for the degeneration
of the stably marked model (cf.\ figure \ref{typ}).
We assume that the cover $C \la \PK$ is given birationally by $Z_0^p=f(X_0)$ and
can choose $f(X_0)$ of the form
$f(X_0)=1+b_2X_0^2+b_3X_0^3+b_4X_0^4+X_0^5 \in R[X_0]$.
We write $Q_8$ for the quaternion group and $D_8$ for the dihedral group of order $8$.

\begin{figure}[htb]  \begin{center}
\input{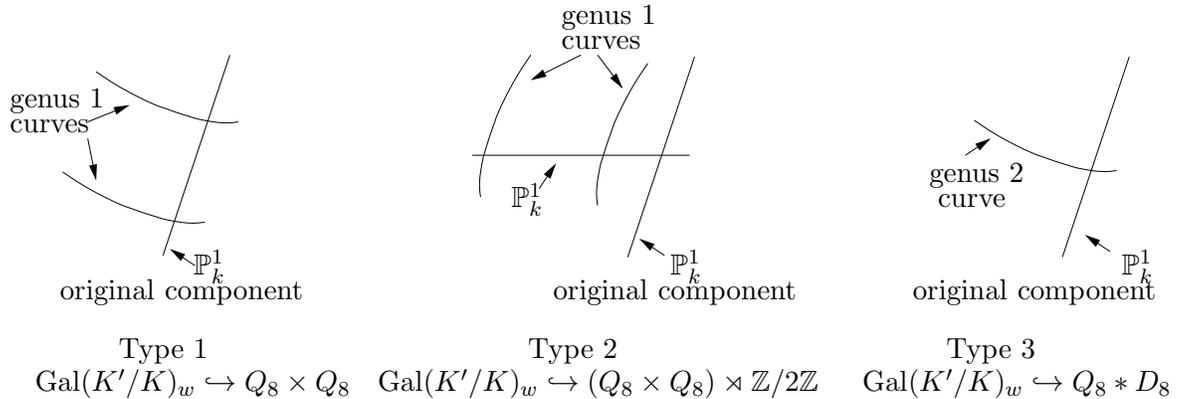}
\caption{The three types of degeneration}
\label{typ}
\end{center}
\end{figure}

We will see that the wild monodromy can be maximal in any of the three types
of degeneration.

\noindent a)
Let $f(X_0)=1+2^{3/5}X_0^2+X_0^3+2^{2/5}X_0^4+X_0^5$ and $K=\Q_2^\ex(2^{1/15})$.
Using our Theorem \ref{algo} and Magma calculations one shows that
$C$ has stably marked reduction of type 1.
Further the finite monodromy
extension is maximal; its Galois group is $Q_8 \times Q_8$.\\

\noindent b) 
Let $K=\Q_2^\ex(a)$ with $a^9=2$. 
Consider the cover $C \la \PK$ given by

$$Z_0^2=f(X_0)=1+a^3X_0^2+a^6X_0^3+X_0^5.$$

\noindent This cover has stably marked reduction of type 2 and maximal finite
monodromy. Its Galois group is isomorphic to
$(Q_8 \times Q_8) \sdp \Z/2\Z$, with $\Z/2\Z$ acting by exchanging the two factors.
\\

\noindent c)
Let $C \la \PK$ be given by the Kummer equation
$Z_0^2=f(X_0)=1+X_0^4+X_0^{5}$ with $K=\Q_2^{{\rm \scriptsize ur}}$. 
This cover has potentially good reduction (i.e.\ is of type 3) and  maximal 
wild monodromy with group the central product $Q_8 * D_8$.
\end{exa}


\begin{remark} \rm
In \cite{Si-Za1} and \cite{Si-Za2} Silverberg and Zarhin study the 
finite monodromy groups for
abelian surfaces. In particular, they classify the finite groups
which can occur as finite monodromy groups for abelian surfaces.
Their examples are mostly in equal characteristic, and they ask at the end
of loc.\ cit.\ section 1 for examples in mixed characteristic. Our Example
\ref{elli} is such a case.
\end{remark}

\medskip

\noindent{\bf Thanks.}
We would like to thank Qing Liu for illuminating discussions. 
Part of this paper was written while the second author was a guest of the 
Max-Planck-Institut in Bonn. We would also like to thank the referees for 
their constructive remarks concerning the presentation of the paper.

\vskip1cm

\begin{flushleft}
Claus LEHR\\
Dipartimento di Matematica Pura ed Applicata\\
Universit\`a degli Studi di Padova,
Via G.Belzoni 7,
35131 Padova, Italia\\
e-mail : {\tt lehr@math.unipd.it}\\
\vskip.5cm
Michel MATIGNON \\
Laboratoire de Th\'eorie des Nombres
et d'Algorithmique Arithm\'etique,
UMR 5465 CNRS \\
Universit\'e de Bordeaux I,
351 cours de la Lib\'eration, 
33405 Talence Cedex, France \\
e-mail : {\tt matignon@math.u-bordeaux1.fr}

\end{flushleft}

\end{document}